\documentclass[11pt]{article} 
\usepackage{amssymb,xypic,latexsym}

\def\neweq{\setcounter{equation}{0}}
\newtheorem{theorem}[equation]{Theorem}
\newtheorem{proposition}[equation]{Proposition}
\newtheorem{corollary}[equation]{Corollary}

\newtheorem{lemma}[equation]{Lemma}

\def\question{\noindent\textbf{Question.} }
\def\remark{\noindent\textbf{Remark.} }
\def\proof{\smallskip\noindent {\it Proof. \ }}

\def\beq{\begin{equation}}

\def\eeq{\end{equation}}

\def\C{\mathbb C}

\def\Z{\mathbb Z }
\def\OO{\mathbb O }

\def\dim{\mbox{dim }}

\def\L{{\cal L}}

\def\LL{{\mathfrak L}}

\def\<{\langle\kern-.08cm\langle}
\def\>{\rangle\kern-.08cm\rangle}

\def\b{{\mathfrak b}}

\def\g{\mathfrak{g}}

\def\m{\mathfrak{m}}

\def\e{{e}}
\def\f{{f}}

\def\IH{\mbox{\sl IH}^{^{\,}{\bullet}}}
\def\IHO{\mbox{\sl IH}^{^{\,}0}}
\def\cg{\C[\g]}
\def\ad{{\mbox{\sl ad}^{\,}}}
\def\rk{{\mbox{\sl rk}^{\,}}}
\def\Hom{{\mbox{\sl Hom}^{\,}}}
\def\Ext{{\mbox{\sl Ext}}}
\def\Specm{{\mbox{\sl Spec}\,}}
\def\ttt{{\mbox{\footnotesize{{\sl {simple}} $G${\sl{-modules}} $V$}}}}
\def\gde{\gd^\e}
\def\gdh{\g^h}
\def\nreg{{\cal N}^{^{{reg}}}}
\newcommand{\ab}{{\hskip 8mm}}
\newcommand{\codim}{\mbox{\sl {codim}}\,}

\newcommand{\rs}{{\mbox{\tiny reg}}}
\newcommand{\gr}{{\mathsf {Gr}}}
\newcommand{\gd}{\check{\mathfrak{g}}}
\newcommand{\hd}{\check{\mathfrak{t}}}
\newcommand{\sll}{\mathfrak{s}\mathfrak{l}_2}
\newcommand{\Lie}{\mbox{\sl Lie}\ }
\newcommand{\ann}{\mbox{\sl Ann}}
\newcommand{\ddim}{\mbox{\sl dim}_{_{\,\C}}}

\newcommand{\ddum}{\mbox{\sl dim}^{\,}}
\def\triv{{\mathbf {triv}}}
\newcommand{\sset}{{\subset}}

\newcommand{\deltav}{{\delta_{_V}}}
\date{}

\title{%
\large{\textbf{\huge{Loop Grassmannian cohomology,
the principal nilpotent and Kostant theorem}}}
\author{\large
VICTOR GINZBURG
}}

\begin{document}

\maketitle

\begin{abstract}
Given a complex projective algebraic variety,
write $H^\bullet(X, \C)$ for its cohomology with complex coefficients
and $\IH(X, \C)$ for its Intersection cohomology.
We first show that under some fairly general conditions 
the canonical map $H^\bullet(X, \C) \to \IH(X, \C)$ is injective.

Now let $\gr := G((z))/G[[z]]$ be the loop Grassmannian for a
complex semisimple group $G$, and let $X$ be the closure
of a $G[[z]]$-orbit in $\gr$.
We prove, using the general result above, a conjecture of D. Peterson
describing the cohomology algebra  $H^\bullet(X, \C)$ in terms of the
centralizer of the principal nilpotent in the Langlands dual of 
$\Lie(G)$.

In the last section we give a new "topological" proof of Kostant's theorem
about the polynomial algebra of a semisimple Lie algebra, based on purity
of the  equivariant intersection cohomology groups of $G[[z]]$-orbits on $\gr$.
\end{abstract}

\section{Main results.}  \label{main}

\ab The purpose of this note is to prove a result relating
the cohomology of some Schubert varieties in the affine Grassmannian
to the centralizer of a
principal nilpotent in the Langlands dual semisimple Lie algebra.
This result was communicated to me, as a conjecture, by Dale 
Peterson in Summer 1997. 

 Our proof of  Peterson's conjecture is based on a general
geometric result about intersection cohomology of algebraic varieties
with a $\C^*$-action, which we explain now.

 Let $X$ be a smooth complex  projective variety with an algebraic 
$\C^*$-action. Assume $\C^*$ acts on $X$
with isolated fixed points, and write $W$ for the fixed point set.
For each fixed point $w\in W$
let
$$ C_w = \{ x\in X\enspace\mid\enspace
\lim_{t \to 0}\, t\cdot x = w\}\quad,\quad
 t\in \C^*,$$ 
denote the corresponding attracting set,
where $t\cdot x$ stands for the action of $t$ on $x\in X$.
These sets form the Bialynicki-Birula cell-decomposition 
$X = \sqcup_w\, C_w$, see [BB]. 

Fix $w\in W$, let $j: C_w \hookrightarrow X$, be the inclusion,
$X_w=\overline{C_w}$ the closure of the cell,
and $IC(X_w,\C)$ the corresponding
intersection cohomology complex.
We use the "naive" normalization in which the restriction
of $IC(X_w,\C)$ to $C_w$ is the constant sheaf concentrated in degree 0 (not in degree
$-\ddim X_w $ as in \cite{BBD}). Thus non-zero cohomology sheaves, ${\cal H}^iIC(X_w,\C)$,
may occur only in degrees $0\leq i\leq \ddim X_w$; in particular, we have
${\cal H}^0IC(X_w,\C) = {\cal H}^0j_*\C_{\,C_w}\,,$
is the direct image of the constant sheaf on $C_w$.
Using  standard truncation functors, $\tau_{_{\leq j}}, \, j\in {\mathbb Z}$, 
on the derived category,
see e.g., \cite{BBD}, we may rewrite this isomorphism
in the form $\tau_{_{\leq 0}}IC(X_w,\C)={\cal H}^0j_*\C_{\,C_w}$.
Therefore,
 one obtains by adjunction canonical morphisms
\begin{equation}\label{deg}
\C_{X_w} \to {\cal H}^0j_*\C_{\,C_w}=
\tau_{_{\leq 0}}IC(X_w,\C) \to IC(X_w,\C)\,.
\end{equation}
The composition of these morphisms induces a natural map on hyper-cohomology
\begin{equation}\label{kappa}
\varkappa:\; H^\bullet(X_w,\C) = H^\bullet(\C_{X_w})\, \longrightarrow\, H^\bullet(IC(X_w,\C))=\IH(X_w, \C)\,.
\end{equation}

Our general geometric
result is

\begin{theorem}  \label{thm1}
Assume that 
the decomposition $X=\sqcup_{w\in W}\, C_w $ is a
stratification of
$X$.  
Then, for any $w\in W$, the map
$\varkappa: H^\bullet(X_w,\C) \rightarrow \IH(X_w,\C)$ is injective, or equivalently,
the dual map to homology
$\IH(X_w,\C) \rightarrow H_\bullet(X_w,\C)$ is surjective.
\end{theorem}

\vskip 3mm

We now recall some basic notation concerning loop groups.
 Let $\C((z))$ be the field of formal
Laurent power series, and $\C[[z]] \sset\C((z))$ its ring of integers,
that is the ring of
formal power series regular at $z=0$. Fix a complex connected
semisimple group $G$ with trivial center,
i.e., of adjoint type, and write $G((z))$, resp. $G[[z]]$,
for the set of its $\C((z))$-rational, resp. $\C[[z]]$-rational,
points. The coset space $\gr := G((z))/G[[z]]$ is called
the {\it loop Grassmannian}. The space $\gr$
has a natural structure of a direct limit of a sequence
of projective varieties of increasing dimension, see
 e.g. \cite[\S 1.2]{G2} or \cite{Lu}. Furthermore,  all orbits
of the left $G[[z]]$-action on $\gr$ are finite dimensional.
Choosing a maximal torus  and a Borel subgroup
$T\sset B\sset G$ gives a natural labelling of
the $G[[z]]$-orbits in $\gr$ by anti-dominant coweights 
$\lambda \in \mbox{Hom}(\C^*,T)$. We write $\OO_\lambda$ for the 
$G[[z]]$-orbit  corresponding to an anti-dominant coweight
$\lambda$. The closure, $\overline{\OO}_\lambda\sset \gr$ is known
to be a finite dimensional projective variety, singular in general.

Let $\gd$ be the complex semisimple Lie algebra dual to
$\Lie G$ in the sense of Langlands. That is, $\gd$ has 
a Cartan subalgebra $\hd$ identified with $(\Lie T)^*$, the
dual of $\Lie T$, and the root system of $(\gd, \hd)$
is dual to that of $(G,T)$. Thus, the coweight lattice 
$\mbox{Hom}(\C^*,T)$ becomes identified canonically 
with the weight lattice $\hd_{_\Z}^* \sset \hd^*$.

\vskip 1pt
Fix a
principal $\sll$-triple $\langle h,\e,\f \rangle \subset \gd$,
such that $h\in \hd$ and such that $\e \in \gd$ is a principal
nilpotent contained in the span of {\it positive} root vectors.
Then the centralizer of $h$ in $\gd$ equals $\hd$.
Further, the centralizer algebra $\gde$ is an
abelian Lie subalgebra in $\gd$
whose dimension equals  $\dim \hd$. The space
$\gde$ is stable under the adjoint $h$-action on $\gd$, and the
weight decomposition with respect to $\mbox{ad}\,h$ puts a grading
on $\gde$. Because
$\gde$ is abelian we identify the enveloping algebra $U(\gde)$ with
the symmetric algebra $S(\gde)$, and view it as a graded algebra 
with the grading induced from that on  $\gde$.
The following
result has been proved in [G2, Proposition 1.7.2] (and
independently proved by Peterson in the simply-laced
case).

\begin{proposition}  \label{H(Gr)}
There is a natural graded algebra isomorphism $\varphi:
H^\bullet(\gr, \C) 
\stackrel{\sim}{\longrightarrow}
U(\gde)$.
\end{proposition}

Given an anti-dominant weight $\lambda \in \hd_{_\Z}^*$,
let $V_{\lambda}$ denote the irreducible
representation of $\gd$ with lowest weight $\lambda$. Choose $v_{\lambda}$,
a lowest weight vector, and
write $\ann_{ U(\gde)}(v_{\lambda})$ for the annihilator
of $v_{\lambda}$ in $U(\gde)$. On the other hand let $i: \overline{\OO}_\lambda
\hookrightarrow \gr$ be the imbedding of
the projective variety $\overline{\OO}_\lambda$ labelled by $\lambda$.

With this understood, 
the result conjectured by Peterson reeds

\begin{theorem} \label{peterson}
The restriction map $i^*: H^\bullet(\gr, \C) \to H^\bullet(\overline{\OO}_\lambda,
\C)$ is surjective and induces, via Proposition \ref{H(Gr)},
a graded algebra isomorphism
$$H^\bullet(\overline{\OO}_\lambda, \C) \simeq U(\gde)/\mbox{\ann}_{ U(\gde)}(v_{\lambda})\;.$$
\end{theorem}

  This generalizes \cite[Proposition 1.9]{G2} as well as the result \cite[Proposition 1.8.1]{G2},
 saying
that if $\lambda$ is minuscule then $U(\gde) \cdot v_{\lambda} =
V_{\lambda}$.  Our geometric proof of Theorem \ref{peterson} also implies the following

\begin{corollary} \label{peterson2}
For any vector $v\in V_\lambda$ one has: $\enspace\mbox{\ann}_{ U(\gde)}(v)
\supseteq \mbox{\ann}_{ U(\gde)}(v_{\lambda})\,.\enspace\square$
\end{corollary}

\section{Proof of Theorem 1.3.}  \label{proof1}\neweq
The strategy of the proof follows the pattern of \cite[\S 3]{G1}.

First of all, we enumerate the strata
$\{C_w\}_{w\in W}$ in a convenient way. To that end, write $\xi$ for the vector field
on $X$ generating the $S^1$-action on $X$ arising from the $\C^*$-action by restriction to the unit
circle. Choose an $S^1$-equivariant  K\"ahler form $\omega$ on $X$, and let $i_\xi\omega$ be the
1-form obtained by contraction. This form is exact since $H^1(X, \C)=0\,,$ hence, there is a function
$f\in C^\infty(X)$ such that $i_\xi\omega=df$. The function $f$ is known to be a Morse function whose critical
points are precisely the fixed points of the $\C^*$-action. Moreover, the Bialynicki-Birula decomposition
coincides with the cell-decomposition associated to $f$ by Morse theory, see e.g. \cite[ch.2]{CG}.
 We enumerate all the fixed points $\{w_1,\ldots,\,w_N\}= W$
in such a way that $f(w_1) \leq f(w_2) \leq\ldots\leq f(w_N),$ and put
$C_n:= C_{w_n},\,n=1,\ldots,N$. 
The sets $X_k:=\sqcup_{n\leq k}\, C_n$ form an increasing filtration of $X$
by closed algebraic subvarieties. Note that $X_{w_k}$ is an irreducible component of $X_k$, so that 
$IC(X_{w_k}, \C)$ is a direct summand of $IC(X_k, \C)$. Therefore, we may (and will) replace
$X_{w_k}$ by $X_k$ in some arguments below.

In addition to the Bialynicki-Birula decomposition $X=\sqcup_{w\in W}\, C_w$ considered so far,
which is often referred to as the {\it plus-decomposition}, one also has the dual
{\it {minus-decomposi}{
tion}} $X = \sqcup_{w\in W}\, C^-_w$, where 
$$
C^-_w = \{ x\in X\enspace\mid\enspace
\lim_{t \to \infty}\, t\cdot x = w\} $$
is the repulsing set at $w\in W$.
Let $\overline{C^-_w}$ denote the closure of $C^-_w$, and write
$c_n\in H^\bullet(X, \C)$ for the Poincar\'e dual of the fundamental class of $\overline{C^-_{w_n}}$.
 Recall (see e.g. \cite[p.488]{G1}) that
 the closure $\overline{C^-_n}$ does not intersect
$X_{n-1}$ and meets $X_n$ transversally in a single point, $w_n$. 
Therefore, $\langle c_n, [C_n]\rangle=1$, and the
 classes $\{c_n\}_{ n=1,\ldots,N}$  form a basis of $H^\bullet(X, \C)$.

From now on we fix some $w=w_k\in W$ and 
let $i_w: X_w \hookrightarrow X$ denote the inclusion. It follows that
the classes $i^*_wc_n$ such that $w_n\in X_w$ form a basis of $H^\bullet(X_w, \C)$.
Abusing the notation we will often write $c_n$ instead of  $i^*_wc_n$.

Proving the theorem amounts to showing that, for  all $n$ such that $w_n\in X_w$, 
the classes in $\varkappa(c_n)\in\IH(X_w,\C)$ are linearly independent.
Assume to the contrary, that there exists a non-trivial linear relation:
\begin{equation}\label{relation}
\sum_{\{n\,|\,w_n\in X_w\}}\,
\lambda_n\cdot\varkappa(c_n) =0\,.
\end{equation}
Let
$n$ be the minimal index such that $\lambda_n\neq 0$. We keep this choice of $n$ from now on.
Put $d=2\ddim X_n=2(\ddim X-\ddim C^-_n) $,
so that the cohomology class $c_n\in H^\bullet(X)$ has degree $d$.

We have  natural diagrams of inclusions
\begin{equation} \label{eq1}
i_n:\,X_n\, \hookrightarrow\, X\quad,\quad
X_{n-1}\, \stackrel{v}{\hookrightarrow}\,X_n\, \stackrel{u}{\hookleftarrow}\,C_n\,.
\end{equation}
Writing 
$\LL:= IC(X_w)$ for the intersection complex on $X_w$ we get canonical morphisms
$\quad$
$\LL \to (i_n)_*i^*_n\LL \to (i_nu)_*(i_n u)^*\LL$. These sheaf morphisms induce natural
 maps
on hyper-cohomology:
$$\IH(X_w) = H^\bullet\LL \to H^{\bullet}(i^*_n\LL)\to H^{\bullet}(u^*i^*_n\LL)\,.$$
Observe that the cohomology class $c_n\in H^d(X)$ acts by multiplication on each of the 
hyper-cohomology groups above. We consider the following diagram, see \cite[(3.8a)]{G1}:
\begin{equation}\label{diagram}
\diagram
H^\bullet(X_w)\rto^{\varkappa}\dto_{c_n\cup} &
H^\bullet\LL \rto^{i^*_n}\dto_{c_n\cup} & H^{\bullet}(i^*_n\LL)\rto^{u^*}\dto_{c_n\cup} &
H^{\bullet}(u^*i^*_n\LL) \ddouble_{c_n\cup}
\\
H^{\bullet+d}(X_w)
\rto^{\varkappa}&
H^{\bullet+d}\LL \rto^{i^*_n} & H^{{\bullet+d}}(i^*_n\LL) & H^{{\bullet+d}}_c(u^*i^*_n\LL)\lto_{u_!}
\enddiagram
\end{equation}
The group $H^{{\bullet+d}}_c(-)$ at the bottom right corner of the diagram
stands for the cohomology with compact support,
the rightmost vertical map is essentially
the standard Thom isomorphism $H^0({\mathbb R}^d) \stackrel{\sim}{\longrightarrow} H^d_c({\mathbb R}^d)$, and
the maps $u^*$ and $u_!$ are induced by the inclusion $u$ in (\ref{eq1}). The first two squares in (\ref{diagram})
clearly commute. Further, for any constructible complex
$\L$ on $X_n$ which is constant along the stratification
$X_n=\sqcup_{j\leq n}\, X_j$, the action on $H^{\bullet}(\L)$ of the Poincar\'e dual of the
fundamental class of the submanifold $\varepsilon: C^-_n\hookrightarrow X$ is given by the composition
of the following natural maps (see, e.g. proof of the `hard Lefschetz theorem' in \cite{BBD}):
$$H^{\bullet}(\L)\,\stackrel{\varepsilon^*}{\longrightarrow}\, H^{\bullet}(\varepsilon^*\L)
 \simeq H^{{\bullet+d}}(\varepsilon^!\L) \,\stackrel{\varepsilon_!}{\longrightarrow}\,
H^{{\bullet+d}}(\L)\,. $$
In the case $\L=i^*_n\LL$ the composition above amounts to going, in diagram (\ref{diagram}),
 along the arrow
$u^*$  followed by the rightmost vertical arrow, and finally along
the arrow $u_!$. This shows that
the right square in (\ref{diagram}) commutes. Thus, (\ref{diagram}) is a commutative diagram.

We will make use of the following result, due to Soergel \cite[Lemma 19]{S}:

\begin{lemma} \label{soergel}
The map $u_!$ in diagram (\ref{diagram}) is injective.$\enspace\square$
\end{lemma}

\remark This result was proved in \cite{S} by showing that
the hyper-cohomology long exact sequence associated to the distinguished
triangle
$u_!\,u^!\L \to \L \to v_*\,v^*\L\,$ splits, provided $\L$ is pointwise pure.
The pointwise purity of $\LL$ (as well as the above Lemma) was verified in
\cite[Lemma 3.5]{G1} under the additional technical condition
\cite[(1.2)]{G1}. This additional condition is however not necessary and can be avoided as follows.
One first argues that, since  $C^-_n$ is an  algebraic subvariety transverse to all the
strata $C_j$, the restriction $\varepsilon^*\LL$ is pure. The result then follows by a standard
argument as, e.g., in the proof of \cite[Lemma 3.5]{G1}. The reason we assumed condition
 \cite[(1.2)]{G1} was that in \cite{G1}, in addition to the injectivity
of the map $u_!$, we also used {\it surjectivity} of the map
$u^*$ in diagram (\ref{diagram}). That surjectivity plays no role in the present paper.
$\enspace\lozenge$
\smallskip

We observe first that $1\in H^0(X_w)$ and 
we have $\varkappa(c_n)=c_n\cup\varkappa(1)$. Using diagram (\ref{diagram}) we find
\begin{equation}\label{eq2}
i^*_n\varkappa(c_n)=c_n\,\cup\, i^*_n\,\varkappa(1)=u_!(c_n\,\cup\,u^*\,
 i^*_n\,\varkappa(1))\,.
\end{equation}
Now, it is immediate from (\ref{deg}) that 
the class $u^*\, i^*_n\,\varkappa(1)$ is a generator of the 0-th hyper-cohomology group
of the complex $u^*i^*_n\LL$, hence non-zero. Therefore, the class
$c_n\,\cup\,u^*\, i^*_n\,\varkappa(1)$ is again non-zero, by the Thom isomorphism.
Hence, the RHS of (\ref{eq2}) is non-zero, due to Lemma \ref{soergel}.
We conclude that 
\begin{equation}\label{eq3}
i^*_n\varkappa(c_n)\neq 0
\end{equation}

To complete the proof of Theorem \ref{thm1} 
we apply the map $i^*_n$ to the linear relation (\ref{relation}). Bearing in mind our choice of $n$
we obtain
$$0=i^*_n\Bigl(\sum_{\{m\,|\,w_m\in X_w\}}\,
\lambda_m\cdot\varkappa(c_m)\Bigr)=
\sum_{m \geq n}\,\lambda_m\cdot i^*_n\,\varkappa(c_m) =
\lambda_n\cdot i^*_n\,\varkappa(c_n)\,,$$
where the last equality is due to the fact that, for any $m>n$, the fundamental
class of $\overline{C^-_m}$ does not intersect $X_n$, whence $i^*_nc_m = 0$.
Thus, the equation above yields $\lambda_n\cdot i^*_n\,\varkappa(c_n)=0$,
and in view of (\ref{eq3}) we deduce $\lambda_n=0$. The contradiction
completes the proof of the Theorem. $\enspace\square$
\medskip

\question We do not know whether Theorem 1.3 is a formal consequence 
of the main theorem of 
[G1], in view of
 the similarity between the proofs  of the two theorems. 
$\enspace\lozenge$

\medskip

It often happens in applications that the
$\C^*$-action on $X$ can be extened to an algebraic action of a complex
torus $T\supset \C^*$. Then each stratum of the Bialynicki-Birula decomposition 
$X = \sqcup_w\, C_w$ is $T$-stable since the actions of $T$ and $\C^*$ commute.
Therefore, for any $w\in W$, we may consider $T$-equivariant cohomology groups
$H^\bullet_T(X_w,\C)$ and $T$-equivariant intersection cohomology groups
$\IH_T(X_w,\C)$, cf. \cite[8.3]{G2}.

\begin{corollary}  
If 
the decomposition $X=\sqcup_{w\in W}\, C_w $ is a
stratification of
$X$ then, for any $w\in W$, the natural map
$H^\bullet_T(X_w,\C) \rightarrow \IH_T(X_w,\C)$ is injective.
\end{corollary}

\proof
It is known that both $H^\bullet_T(X_w,\C)$ and $\IH_T(X_w,\C)$ are finitely
generated modules over $H^\bullet_T(pt)\simeq \C[\Lie T]$. Hence to prove injectivity
it suffices to show that, for any maximal ideal ${\mathfrak m}\subset \C[\Lie T]$,
the localized map $H^\bullet_T(X_w,\C)_{({\mathfrak m})} \rightarrow \IH_T(X_w,\C)_{({\mathfrak m})}$ is injective.
Any maximal ideal in $\C[\Lie T]$ consists of the polynomials vanishing at a given point
$t\in \Lie T$. Therefore, we must show that, for any $t\in \Lie T$, the
localized map $H^\bullet_T(X_w,\C)_t \rightarrow \IH_T(X_w,\C)_t$ is injective.
But the latter map may be replaced, due to the {\it Localization theorem} in equivariant cohomology,
cf., \cite[Thm.8.6]{G2},
by a similar map,
$H^\bullet(X_w^t,\C) \rightarrow \IH(X_w^t,\C)$, between the corresponding
non-equivariant cohomology groups of the
$t$-fixed point set, $X_w^t$. The result now follows from Theorem \ref{thm1},
applied to the $\C^*$-manifold $X^t$.$\enspace\square$
\vskip 2mm

\section{The loop Grassmannian.}
 \label{proof2}\neweq

We would like to apply Theorem \ref{thm1} to $X=\gr$, the loop Grassmannian.

Recall that we have fixed $T\subset B \subset G$, a maximal torus and a Borel subgroup in $G$.
Define an Iwahori subgroup $I \subset G[[z]]$ to be formed by all loops
$f\in G[[z]]$ such that $f(0)\in B$. It is known, see e.g., \cite{Lu}, that  $I$-orbits
form a cell-decomposition of $\gr$ that refines the stratification by $G[[z]]$-orbits,
$\gr=\sqcup_\lambda\,\OO_\lambda$. In particular, for any $\lambda$, the variety
$\overline{\OO}_\lambda$ is the closure of a single $I$-orbit. It is known further that the
decomposition of $\gr$ into $I$-orbits coincides with the Bialynicki-Birula decomposition,
$\gr =\sqcup_\lambda\,C_\lambda\,,$
with respect to an appropriate one-parameter subgroup $\C^*\subset T$. Thus, we are in the
setup of Theorem \ref{thm1} except that the variety $\gr$ is neither finite-dimensional,
nor smooth.

There is a standard way, see e.g., \cite{KT}, to go around this difficulty. Specifically,
the space $\gr$ may be imbedded into a slightly larger infinite dimensional variety
${\mathsf {\widetilde{Gr}}}$, which is a union of $G[z^{-1}]$-orbits  of finite {\it codimension}.
Thus, ${\mathsf {\widetilde{Gr}}}$ has the structure of a direct limit of infinite-dimensional
smooth open subsets, hence may be regarded as a {\it smooth} variety
(see \cite{KT} or \cite[\S\S 6.1-6.4]{G2} for more details about such smooth infinite dimensional
varieties). Although the variety ${\mathsf {\widetilde{Gr}}}$ is by no means
compact, there is an explicit {\it minus-decomposition}
${\mathsf {\widetilde{Gr}}}
=\sqcup_\lambda\,C^-_\lambda$, cf., \cite[6.4]{G2}, that enjoys all the properties of the minus decompostion
for a $\C^*$-action on a smooth projective variety,
 that were exploited in the proof
of Theorem \ref{thm1} above. Therefore the proof of the theorem goes through.
We conclude that  theorem  \ref{thm1} holds for $\gr$ in the sense that,
for any $\lambda$, the canonical map
$H^\bullet(\overline{\OO}_\lambda,\C) \rightarrow \IH(\overline{\OO}_\lambda,\C)$ is injective.

We now recall the main result of \cite{G2}. 
Let $P(\gr)$ be the category 
of semisimple $G[[z]]$-equivariant 
perverse sheaves on $\gr$ with compact support.
Also, write $Rep (G^\lor)$ for
the category of finite dimensional representations of  $G^\lor$, the Langlads dual of $G$.
Then we have (see \cite[Theorem 1.4.1 and Theorem 1.7.6]{G2}):

\begin{theorem}\label{tensor} (i)
There is an equivalence of the categories $P(\gr)$
and  $Rep(G^\lor)$
which sends $ IC(\overline{\OO}_\lambda, \C)$ to $V_\lambda$. 

(ii) For any $\L \in P(\gr)$, the hyper-cohomology $H^\bullet(\L)$ gets identified,
under the equivalence, with the underlying vector space of the corresponding representation
of $G^\lor$.

(iii) Furthermore, for any
$u\in H^\bullet(\gr, \C)$,
the natural action of $u$ on the hyper-cohomology $H^\bullet(\L)$ corresponds,
via (ii) and the isomorphism $\varphi:
H^\bullet(\gr, \C) 
\stackrel{\sim}{\longrightarrow}
U(\gde)$ of Proposition \ref{H(Gr)},
to the natural
action of $\varphi(u)\in U(\gde)$ in the corresponding $G^\vee$-module.
\hfill$\enspace\square$\break
\end{theorem}

\noindent
{\textbf {Proof of Theorem 1.5.}} Fix $\lambda$, and let
$\IH(\overline{\OO}_\lambda, \C)\stackrel{\sim}{\longrightarrow}
V_\lambda$ be the identification of Theorem \ref{tensor}(ii).
 According to \cite{G2} this map  sends the unit,
$1\in \IHO(\overline{\OO}_\lambda, \C)$,
to a lowest weight vector $v_\lambda\in
V_\lambda$. Therefore, Theorem \ref{tensor} implies that the map
$\varphi$ of \ref{tensor}(iii) induces a graded algebra isomorphism:
\begin{equation}\label{ann}
\varphi:\,
H^\bullet(\gr, \C)/\ann_{H^\bullet(\gr, \C)}(\varkappa(1))
\,\simeq\,U(\gde)/\ann_{ U(\gde)}(v_{\lambda})\,.
\end{equation}

Let $i: \overline{\OO}_\lambda
\hookrightarrow \gr$  denote the imbedding.
Note that the $H^\bullet(\gr , \C)$-action
on $\IH(\overline{\OO}_\lambda, \C)$ factors through the restriction map
$i^*: H^\bullet(\gr, \C)\to H^\bullet(\overline{\OO}_\lambda, \C)$. The
restriction map is surjective since the dual
map on homology $i_*: H_\bullet(\overline{\OO}_\lambda, \C)
\to H_\bullet(\gr, \C)$ is injective (because both spaces have compatible cell-decompositions
by $I$-orbits of {\it even real} dimension). Combining the surjectivity observation above
with 
(\ref{ann}) we obtain the following chain of algebra isomorphisms:
$$
H^\bullet(\overline{\OO}_\lambda, \C)/
\ann_{H^\bullet(\overline{\OO}_\lambda, \C)}(\varkappa(1))\,\simeq\,
H^\bullet(\gr, \C)/\ann_{H^\bullet(\gr, \C)}(\varkappa(1))
\,\simeq\,U(\gde)/\ann_{ U(\gde)}(v_{\lambda})\,.
$$
But $\ann_{H^\bullet(\overline{\OO}_\lambda, \C)}(\varkappa(1))=0$ because of
Theorem \ref{thm1} applied to the variety
  $\gr$. The isomorphism of Theorem 1.5 follows.
$\enspace\square$

\section{"Topological" proof of Kostant's theorem.}
 \neweq
Given a complex connected semisimple group $G$ with Lie algebra $\g$ write
$\cg^G\subset \cg$ for the subring of $\ad G$-invariant polynomials on $\g$.
In \cite{Ko}, B. Kostant established the following fundamental result

\begin{theorem}\label{kostant}
There is a graded $G$-stable subspace ${\mathsf H}\subset \cg$ such that the multiplication 
in $\cg$ gives rise to a vector space isomorphism
$$\cg^G\otimes_{_\C} {\mathsf H}\,\stackrel{{mult\atop \sim}}{\longrightarrow}\,\cg\,.$$
\end{theorem}

We are going to show that this theorem may be viewed as a manifestation of "purity"
for the equivariant intersection cohomology, $\IH_T(\overline{\OO}_\lambda, \C)$.
Recall first, that an element $x\in\g$ is called {\it regular} if $\g^x$, 
the centralizer of 
$x$ in $\g$, has the minimal possible dimension, $\rk\g$. Given $x\in \g$ and a
finite dimensional rational
$G$-module $V$, write $V^{\g^x}$ for the subspace in $V$ annihilated by the subalgebra
$\g^x$.

\begin{proposition}\label{key}
For any finite dimensional rational
$G$-module $V$ whose weights are contained in
 the root lattice of $G$, the function:
 $x\mapsto \ddim V^{\g^x}$ 
is constant on the set of regular elements of $\g$.
\end{proposition}

\remark Note that for $V=\g$, the adjoint representation, the Proposition amounts
to the definition of a regular point.$\enspace\lozenge$
\smallskip

Of course,   Proposition \ref{key}
follows from Theorem \ref{kostant},
by a  well-known argument due to Kostant \cite{Ko}. Our main observation is
that the results of \cite{G1} and \cite{G2} combined together
yield an alternative "topological" proof of  Proposition \ref{key}, independent
of Theorem \ref{kostant}.
\medskip

\noindent
{\textbf {Proof of Proposition \ref{key}:}} The natural projection
$\pi: \g= \Specm\cg \twoheadrightarrow \Specm \cg^{G}$ sets up a bijection
between regular adgoint $G$-orbits in $\g$ and (closed) points of the scheme
$ \Specm \cg^{G} \simeq {\mathfrak t}/W$. Fix a representation $V$, as in the proposition.
The function $\deltav: x\mapsto \ddim V^{\g^x}$
is clearly constant on each $G$-orbit, hence, when restricted to regular elements, it
may (and will) be viewed as a function on $ \Specm \cg^{G}$.
By semicontinuity, the value
of this function at a special (regular) orbit can not be less than its value at the
generic orbit. Observe further that there is a $\C^*$-action on $\g$
by homotheties, preserving the set of regular elements. 
It induces a natural $\C^*$-action on $\Specm \cg^G$
with the origin, $\m_\circ$, being the unique attracting fixed point.
Thus the point $\m_\circ$ is the "most special" point in $\Specm \cg^G$
in the sense that, for any other point $\m\in\Specm \cg^G$, we have
$\deltav(\m)\leq \deltav(\m_\circ)$.
Thus,
it suffices to show that the
value of the function $\deltav$ at $\m_\circ$ equals its generic value.
Note that,
since the centralizer of a generic element is a Cartan subalgebra,
 the generic value of $\deltav$ is equal to $\ddum V(0)$,
the zero-weight multiplicity in $V$.
Thus, we must prove that, if $x$ is a regular nilpotent, then $\ddum V^{\g^x}=\ddum V(0)$.

To this end, we may replace in all the arguments the Lie algebra $\g$ by $\gd$, its
Langlands dual. Thus we let $V_\lambda$ be a simple finite-dimensional $G^\vee$-module
with lowest weight $\lambda$, let $V_\triv$ be
the trivial $G^\vee$-module , and write $\e$ for a regular
nilpotent in $\gd$. We have $\ddum V_\lambda^{\gde}= \ddum\Hom_{\gde}(V_\triv, V_\lambda)$.

The space $\Hom_{\gde}(V_\triv, V_\lambda)$
 may be expressed in terms of the geometry of the loop Grassmannian $\gr$.
 Specifically, Theorem 1.10.3 of \cite{G2} (whose proof depends on
\cite{G1} in an essential way) gives an isomorphism of vector spaces
\begin{equation}\label{ext}
\Ext^\bullet_{D^b(\gr\,)}(IC_\triv, IC(\overline{\OO}_\lambda))\,\simeq\,
\Hom_{\gde}(V_\triv, V_\lambda)\,,
\end{equation}
where $IC_\triv$ is the skyscrapper sheaf on the one-point orbit 
$\,i_\triv: \OO_\triv\hookrightarrow \gr\,$ 
corresponding
to the trivial representation. The LHS of (\ref{ext}) equals, by adjunction, 
$i_\triv^!IC(\overline{\OO}_\lambda)$. Thus, proving the theorem amounts to showing that,
for any anti-dominant $\lambda$ in the root lattice, one has
\begin{equation}\label{e=h}
\ddum H^\bullet i_\triv^!IC(\overline{\OO}_\lambda)\,=\,
\ddum V_\lambda(0)
\end{equation}

We interpret the last equation in terms of equivariant cohomology as follows.
Let $T\subset G$ be the maximal torus whose fixed points in $V$ form the subspace $V(0)$.
The torus $T$ acts naturally on $\gr$ preserving all the strata $\OO_\lambda$.
It follows that both $IC(\overline{\OO}_\lambda)$ and $i_\triv^!IC(\overline{\OO}_\lambda)$
are $T$-equivariant complexes, see e.g., \cite[8.3]{G2}. We may therefore consider
the $T$-equivariant hyper-cohomology,
$H^\bullet_T(i_\triv^!IC(\overline{\OO}_\lambda))$, which is a module over
$H^\bullet_T(pt)=\C[{\mathfrak t}]$.
But the complex $i_\triv^!IC(\overline{\OO}_\lambda)$ is pure, by \cite[Lemma  3.5]{G1}.
Hence the $\C[{\mathfrak t}]$-module $H^\bullet_T(i_\triv^!IC(\overline{\OO}_\lambda))$
is free, by \cite[Theorem 8.4.1]{G2}. Moreover, the geometric fiber of this free module
at the origin $0\in {\mathfrak t}$ is isomorphic to 
$H^\bullet(i_\triv^!IC(\overline{\OO}_\lambda))$, the non-equivariant cohomology,
by \cite[Corollary 8.4.2]{G2}. On the other hand, by the geometric construction of
a fiber functor on $P(\gr)$ given in
\cite[3.9-3.10]{G2}, the fixed point decomposition \cite[3.6]{G2}
on the equivariant intersection
cohomology corresponds, via Theorem \ref{tensor}, to
 the weight decomposition on $V_\lambda$.
Thus, the geometric fiber of $H^\bullet_T(i_\triv^!IC(\overline{\OO}_\lambda))$
at a general point in ${\mathfrak t}$ is precisely the zero-weight subspace, $V(0)$.
Since all fibers of a free module have the same dimension, we conclude
that $\ddum H^\bullet i_\triv^!IC(\overline{\OO}_\lambda)\,=\,
\ddum V_\lambda(0)$, and Proposition \ref{key} follows.$\enspace\square$

\medskip
\noindent
{\textbf {Proposition \ref{key} implies Theorem \ref{kostant}:}}
Given a simple finite dimensional rational $G$-module $V$, let
$\cg^V= \Hom_G(V, \cg)$ denote the $V$-isotypic component of $\cg$. The $G$-action
on $\cg$ being locally finite, one has a $G$-stable direct sum decomposition:
$$\cg = \bigoplus_\ttt\;
V\otimes_{_\C} \cg^V\,.$$
Clearly, for $V=V_{{\mathbf {triv}}}$, we have
$\cg^{V_{{\mathbf {triv}}}}=\cg^G$, and for an arbitrary $V$,
  $\cg^V$ is a graded $\cg^G$-module. 
By the direct sum decomposition above, proving Theorem \ref{kostant}
amounts to showing that, for any $V$, the $\cg^G$-module $\cg^V$ is free.

Thus we may fix a simple $G$-module $V$ whose weights belong to the 
root lattice (otherwise $V$ does not occur in $\cg$),  and concentrate our attention on the 
$\cg^G$-module $\cg^V$. The latter is finitely generated,
by Hilbert's classical result on finite generation of
$G$-invariants, see e.g., \cite{We}. Further, $\cg^V$
is clearly a graded $\cg^G$-module. 
But a finitely generated graded
$\cg^G$-module is free if and only if it is projective. 
To show $\cg^V$ is projective we argue as follows.

View $\cg^V$ as a coherent sheaf on
$\Specm \cg^G$.
Let $\cg^V\!/\m\!\cdot\!\cg^V$ be the geometric fiber of this sheaf
at a closed point $\m\in\Specm \cg^G$, regarded as a maximal ideal in
$\cg^G$. It is known that
$\cg^V$ is a projective $\cg^G$-module if and only if the
function $d_{_V}:\m \mapsto \ddim (\cg^V\!/\m\!\cdot\!\cg^V)$ is constant on the
set of closed points of $\Specm \cg^G$. 
It suffices to show, due to a semi-continuity
argument similar to the one used in the proof of Proposition \ref{key}, that the
value of the function $d_{_V}$ at $\m_\circ$, the origin of  $\Specm \cg^G$,
equals its generic value.

To this end, consider the natural projection
$\pi: \g= \Specm\cg \twoheadrightarrow \Specm \cg^G$. 
If  $\m\in\Specm \cg^G$ is in general position, then $\pi^{-1}(\m)$ is the
single adjoint $G$-orbit through a
semisimple regular element $h\in\g$. This orbit is isomorphic, as a $G$-variety,
to $G/G^h$, where $G^h$ denotes the centralizer of $h$ in $G$. Therefore, a standard
argument involving Frobenius reciprocity, see \cite{Ko} or \cite[\S 6.7]{CG},
 yields
\begin{equation}\label{generic}
\ddim\, (\cg^V\!/\m\!\cdot\!\cg^V) = \ddim\, \C[G/G^h]^V = \ddim\, V^{\gdh}\,.
\end{equation}

On the other hand, we have set-theoretically:
$\pi^{-1}(\m_\circ)= ${\it nilpotent variety of} $\g$,
see e.g. \cite[ch.3]{CG}. The nilpotent variety contains a unique open dense
$G$-orbit, $\nreg$, formed by regular nilpotents. Moreover, the scheme
$\pi^{-1}(\m_\circ)$ is reduced at any point of $\nreg$, see Lemma
\ref{add} below and \cite{Ko}. Hence, 
the scheme imbedding $\nreg\hookrightarrow\pi^{-1}(\m_\circ)$
induces an injection: 
$\cg^V\!/\m_\circ\!\cdot\!\cg^V \hookrightarrow \C[\nreg]^V$. 
Thus, we have $\ddum (LHS) \leq \ddum (RHS)$.
Choose a regular nilpotent $\e\in\nreg$. Then $\nreg \simeq G/G^\e$, and
the Frobenius reciprocity argument
mentioned above yields:
\[\ddum\C[\nreg]^V= \ddum \C[G/G^\e]^V 
= \ddum V^{\g^\e}\,.
\]
Combining this formula with equation (\ref{generic}) and using Proposition
\ref{key}, we obtain
\[
\ddum(\cg^V\!/\m_\circ\!\cdot\!\cg^V)\,\leq
\,\ddum\C[\nreg]^V 
= \ddum V^{\g^\e}=\ddum V^{\gdh}\,.
\]
Thus the function $d_{_V}$ is constant, and Theorem \ref{kostant}
is proved.$\enspace\square$
\medskip

In the course of our proof of Theorem \ref{kostant} we have used
the following result, due to Kostant.

\begin{lemma}\label{add}
The zero fiber, $\pi^{-1}(\m_\circ)$, is reduced at any point of
 $\nreg$.
\end{lemma}

Kostant proved this result by showing
that the generators of $\cg^G$ have linearly independent
differentials
at any point of $\nreg$. The latter has been verified in \cite{Ko2}
by a direct
computation (see \cite[\S6.7]{CG} for a slightly different argument).
We give an alternative proof of the Lemma, inspired by \cite{BL},
which involves no computation and is independent of \cite{Ko2}.
\smallskip

\noindent
{\textbf {Proof of Lemma:}} Set $\widetilde{\g}= \{(x, \b)\;|\;
x\in \b\,,\, \b= \mbox{{\it Borel subalgebra in }}\g\}\,.$ We have the following
commutative diagram, see \cite[\S3.2]{CG}:
\begin{equation}\label{grothendieck}
\diagram
 &\widetilde{\g}\dlto_{\mu}\drto^{\nu}&\\
\g\drto^{\rho} &&{\mathfrak t}\dlto_{\pi}\\
&{\mathfrak t}/W&\enddiagram
\end{equation}
In this diagram, the map $\mu$ is proper,
 $\nu$ is a smooth morphism, and
$\pi$ is a finite flat morphism, since $\C[{\mathfrak t}]$ is free
over $\C[{\mathfrak t}]^W$.

Observe that diagram
(\ref{grothendieck}) induces a morphism
$\psi: \widetilde{\g}\to \g\times_{_{{\mathfrak t}/W}} {\mathfrak t}$. 
Let $\g^{\rs}\subset \g$ be the Zariski open subset of regular (not necessarily
semisimple)
elements, and $\widetilde{\g}^{\rs} := \mu^{-1}({\g}^{\rs})$. We claim that
the morphism $\psi$ gives an isomorphism:
\begin{equation}\label{eq5}
\psi: \widetilde{\g}^{\rs} \stackrel{\sim}{\longrightarrow}
 \g^{\rs}\times_{_{{\mathfrak t}/W}} {\mathfrak t}\,.
\end{equation}

To prove the claim, note first that since ${\mathfrak t}$ 
is finite and flat over ${\mathfrak t}/W$, the scheme
 $\g^{\rs}\times_{_{{\mathfrak t}/W}} {\mathfrak t}$
is finite and flat over ${\g}^{\rs}$. 
Further, the map $\,\mu:\widetilde{\g}^{\rs} \to {\g}^{\rs}\,$ is proper
and has finite fibers, hence this is a finite morphism. Moreover,
being  a
dominant morphism between smooth schemes of the same dimension,
this  morphism is flat.
Thus, both the sourse and target schemes in (\ref{eq5})
are finite  flat schemes over $\g^{\rs}$, a smooth variety. Therefore,
both schemes are Cohen-Macaulay (see \cite{BL} or \cite[\S2.2]{CG}),
hence, to show that
the map $\psi$ in (\ref{eq5}) is an isomorphism, it suffices to verify that it is
an isomorphism outside
a codimension two subvariety. Let $\g'\subset\g$ be the set of all elements
whose semisimple part is either regular, or belongs to at most one
root hyperplane in a Cartan subalgebra. Then $\codim(\g\setminus\g')
\geq 2$. On the other hand, proving isomorphism (\ref{eq5})
for $\g'$ amounts, effectively,
to an ${\mathfrak s}{\mathfrak l}_2$-computation, which is left to the reader.
This proves (\ref{eq5}).

Now, using (\ref{eq5}), we may identify the smooth morphism
$\nu: \widetilde{\g}^{\rs} \to {\mathfrak t}$ with the projection
$\g^{\rs}\times_{_{{\mathfrak t}/W}} {\mathfrak t}\to$
${\mathfrak t}$. Applying to this projection 
the base change with respect to the flat
map ${\mathfrak t}\to{\mathfrak t}/W$ in diagram
(\ref{grothendieck}), we deduce that the morphism
$\rho: \g^{\rs} \to {\mathfrak t}/W$ is also smooth.
Hence its zero-fiber, $\nreg$, is reduced, and the Lemma follows.
$\enspace\square$

\vskip 4mm
\noindent
\footnotesize{
{\sc Acknowledgements.} {\it  I am grateful to D. Peterson and 
B. Kostant for interesting discussions.}

}
\vskip 1cm
\noindent
\footnotesize{\hphantom{x}\ab
 Department of Mathematics, University of Chicago, 
Chicago IL
60637, USA; 

\noindent
\hphantom{x}\ab ${\mathbf {ginzburg@math.uchicago.edu}}$

\end{document}

\\
Title: Loop Grassmannian cohomology, the principal nilpotent and Kostant theorem
Author: Victor Ginzburg
Subj-Class: Algebraic Geometry; Representation Theory; Differential Geometry
Comments: 10 pages, LaTeX2e
Notes: lemma 4.6 added, title changed, and several misprints corrected.
\\
Given a complex projective algebraic variety, write H(X) for its cohomology with 
complex coefficients and IH(X) for its Intersection cohomology. We first show that, 
under some fairly general conditions, the canonical map H(X)\to IH(X) is injective.

Now let Gr = G((z))/G[[z]] be the loop Grassmannian for a complex semisimple group 
G, and let X be the closure of a G[[z]]-orbit in Gr. We prove, using the general 
result above, a conjecture of D. Peterson describing the cohomology algebra H(X) in 
terms of the centralizer of the principal nilpotent in the Langlands dual of Lie(G).

In the last section we give a new "topological" proof of Kostant's theorem about the 
polynomial algebra of a semisimple Lie algebra, based on purity of the  equivariant 
intersection cohomology groups of G[[z]]-orbits on Gr.
\\

math.AG/9803141
q9dc5
reederma@bc.edu
alg-geom/9711022
ioanid@math.mit.edu